\documentclass[11pt]{amsart}
\usepackage[T1]{fontenc}
\usepackage{lmodern}
\usepackage[margin=1in]{geometry}
\usepackage{amsmath,amssymb,amsthm,mathtools,booktabs,array}
\usepackage[hidelinks]{hyperref}
\numberwithin{equation}{section}
\newtheorem{theorem}{Theorem}[section]
\newtheorem{proposition}[theorem]{Proposition}
\newtheorem{lemma}[theorem]{Lemma}
\newtheorem{corollary}[theorem]{Corollary}
\theoremstyle{definition}
\newtheorem{definition}[theorem]{Definition}
\newtheorem{example}[theorem]{Example}
\theoremstyle{remark}

\DeclareMathOperator{\Gal}{Gal}
\DeclareMathOperator{\Cl}{Cl}
\DeclareMathOperator{\expn}{exp}
\DeclareMathOperator{\Ver}{Ver}
\newcommand{\Z}{\mathbb Z}
\newcommand{\Q}{\mathbb Q}
\newcommand{\ab}{\mathrm{ab}}
\newcommand{\oo}{\infty}
\newcommand{\cL}{\mathcal L}
\newcommand{\cH}{\mathcal H}
\newcommand{\cM}{\mathcal M}

\title{Canonical Quotients and Recovery of Metacyclic $2$-Class Towers over Cyclotomic $\Z_2$-Extensions}
\author{Marina Palaisti}

\begin{document} 	\maketitle
	
	\begin{abstract}
		We study how abelian unramified extensions recover the degree of a finite nonabelian, nonmodular metacyclic $2$-class tower whose Galois group has abelianization $C_2\times C_{2^n}$, $n\ge2$. The quadratic layer determines a canonical quotient, realized by the maximal unramified abelian $2$-extension of one distinguished quadratic overfield. For every $2\le j\le n$, the three intermediate fields of degree $2^j$ in the first abelian layer have the same maximal unramified abelian $2$-extension. These common fields form a canonical filtration reaching the full tower. The corresponding class-group orders satisfy a plateau formula: a single quartic order detects whether the initial quotient is already the full tower, while higher orders recover its exact degree without choosing a chain of intermediate fields. We prove that the resulting bounds are sharp by constructing metacyclic groups with identical complete first-layer transfer diagrams attaining every permitted order. For five explicit real quadratic fields in a family of Chems-Eddin and El Mamry, class-number and capitulation computations determine the full tower groups, with orders $32$, $64$, and $128$, and identify each cyclotomic tower with the base change of its finite-base tower. Throughout, class groups are ordinary and real places are required to remain real.
		
	\end{abstract}

\noindent \textbf{MSC 2020 Subject Class:} 11R23, 11R29, 11R37, 20D15\\
\noindent \textbf{Keywords:} $2$-class field tower, cyclotomic $\Z_2$-extension, Iwasawa module, metacyclic group, Artin transfer, capitulation
	\section{Introduction}
	
	Class field theory identifies the ideal class group of a number field with the Galois group of its maximal unramified abelian extension. Repeating this construction produces the class-field tower, whose Galois group records nonabelian information invisible in the original class group. The tower can be infinite even though the initial class group is finite, as the work of Golod--Shafarevich shows \cite{GolodShafarevich1964}. A natural way to obtain more information is to examine the class groups of unramified overfields and the ideal classes that become principal there. Artin's interpretation of class extension by group transfer \cite{Artin1929} makes these arithmetic observations accessible to finite-group methods.
	
	The use of such methods in $2$-class towers goes back to Taussky \cite{Taussky1937}. For small $2$-class groups, capitulation distinguishes nonabelian possibilities that have the same abelianization; Kisilevsky's treatment of the four-group case \cite{Kisilevsky1976} is a classical example. Benjamin \cite{Benjamin1999} studied second Hilbert $2$-class fields for real quadratic fields with class group $C_2\times C_{2^n}$, $n\ge2$. This is the abelianization considered here. Our question is how much of a metacyclic tower its overfield class groups recover, and which additional fields reveal the information still missing.
	
	We use ordinary class groups. Accordingly, an unramified extension is required to be unramified at finite places and at real places; over a totally real field, all extensions in the tower remain totally real. For a number field or a totally real algebraic extension $B/\Q$, write $\cL(B)$ for its maximal unramified pro-$2$ extension, $\cH(B)$ for its maximal unramified abelian pro-$2$ extension, and
	\[ G(B)=\Gal(\cL(B)/B),\qquad X(B)=\Gal(\cH(B)/B). \]
	For a number field $B$, $X(B)$ is its ordinary $2$-class group. If $B=k_\oo$ is a cyclotomic $\Z_2$-extension, then
	\begin{equation}\label{eq:iwasawa} X(k_\oo)\cong\varprojlim_\nu\Cl_2(k_\nu), \end{equation}
	with norm transition maps; see \cite[Chapter~13]{Washington}. Here $k_\oo/k$ is the cyclotomic extension with Galois group $\Z_2$, and $k_\nu/k$ is its subextension of degree $2^\nu$.
	
	There are two different extension processes in this problem. Moving from $k$ through the fields $k_\nu$ studies class-group growth in the cyclotomic direction. Moving from the fixed field $B=k_\oo$ into its unramified tower studies the nonabelian group $G(B)$. Iwasawa's growth formula \cite{Iwasawa1959} and Greenberg's conjecture \cite{Greenberg1976} concern the first process; our recovery problem concerns the second. In particular, the degree index $j$ used below refers to unramified extensions of $B$, not to cyclotomic levels. Section~\ref{sec:background} explains the maps and finiteness inputs that connect the two processes.
	
	Suppose that $G=G(B)$ is finite and metacyclic, nonabelian and nonmodular, with
	\[ G^{\ab}\cong C_2\times C_{2^n},\qquad n\ge2. \]
	Here $C_d$ denotes a cyclic group of order $d$, and the modular family is defined in Section~\ref{sec:metacyclic}. A metacyclic group has cyclic derived subgroup, so in this setting the unramified tower stops at its second abelian stage; its degree can nevertheless vary. There are three quadratic fields between $B$ and $\cH(B)$. One is distinguished: its corresponding subgroup has noncyclic image in $G^{\ab}$. Denote this field by $K_B$, and write $|X(K_B)|=2^{n+t}$, where $t\ge1$. The extension $\cH(K_B)/B$ need not be the full tower. It is nevertheless canonical: it is the smallest natural quotient extension retaining all three quadratic transfer targets and maps.
	
	The higher layers have a stronger uniformity. For every $2\le j\le n$, there are exactly three fields of degree $2^j$ over $B$ inside $\cH(B)$, and their maximal unramified abelian $2$-extensions coincide. The following theorem states the resulting recovery principle without choosing generators or a chain of intermediate fields.
	
	\begin{theorem}\label{thm:main} Let $B$, $G$, $n$, $K_B$ and $t$ be as above. Put $\cM_1=\cH(K_B)$. For $2\le j\le n$, let $\cM_j=\cH(E)$, where $E$ is any field with \[ B\subseteq E\subseteq\cH(B),\qquad [E:B]=2^j. \] Then $\cM_j$ is independent of the choice of $E$. These fields are Galois over $B$ and form a canonical filtration \[ \cM_1\subseteq\cM_2\subseteq\cdots\subseteq\cM_n=\cL(B). \] There is an integer $k$, $1\le k\le n$, such that \begin{align} [\cM_j:B]&=2^{n+t+\min\{j,k\}},&&1\le j\le n,\label{eq:mainfields}\\ |X(E)|&=2^{n+t-\max\{j-k,0\}},&&2\le j\le n.\label{eq:mainorders} \end{align}
		In particular, for any quartic $E/B$ in $\cH(B)$, \begin{equation}\label{eq:quarticintro} |X(E)|= \begin{cases} |X(K_B)|/2,&\cL(B)=\cM_1,\\ |X(K_B)|,&\cL(B)\ne\cM_1. \end{cases} \end{equation} If $w_1=|X(K_B)|$ and $w_j=|X(E_j)|$ for arbitrary choices of such fields at degrees $2^j$, then 
		\begin{equation}\label{eq:mainrecovery}
			k=\max\{1\le j\le n:w_j=w_1\},\qquad |G|=2^{n+t+k}=2^n w_n. \end{equation} 
	\end{theorem}
	
	The general recovery assertion concerns the canonical filtration and the total tower degree. Determining the group up to isomorphism requires additional recognition hypotheses, as in Corollary~\ref{cor:endpoint}.
	
	The proof is given in Section~\ref{sec:higher}. The order sequence stays constant for $k$ steps, then halves at each step. At the field level the behavior is reversed: the common abelian extension grows by degree two until it reaches the full tower. Both statements are intrinsic. The fields $E_j$ used to read off their orders need not belong to a selected chain.

	\subsection{Relation to earlier work}
	The established transfer--capitulation correspondence is explained in Miyake's survey \cite{Miyake1989}, while Gruenberg--Weiss \cite{GruenbergWeiss2000} study the transfer kernels themselves. In the language of Artin patterns, the observed information consists of transfer targets and kernels. Mayer \cite{Mayer2016} develops their behavior under quotient maps and on descendant trees. This viewpoint suggests the first step of our approach: factor out exactly the subgroup that is invisible to all targets in a specified layer. The quotient construction uses transfer naturality, which we recall explicitly.
	
	The required finite-group structure also has a substantial history. King's canonical presentations \cite{King1973} and the classification of metacyclic $2$-groups by Xu--Zhang \cite{XuZhang2006} explain why relations and admissible parameters, rather than just group orders, matter. For the particular abelianization $C_2\times C_{2^n}$, Aaboun--Zekhnini \cite{AabounZekhnini2023} give subgroup abelianizations and transfer formulas. Their Corollary~3.3 contains the index-two and index-four derived subgroups used here, and their Lemma~3.4 gives the corresponding targets. Section~\ref{sec:metacyclic} states the presentations and reproduces the commutator calculation, so these formulas can be used without consulting their tables. The first-layer restriction and quartic order test are consequences of that existing calculation. The all-degree coincidence of class fields and the canonical recovery filtration are the subsequent conclusions developed here.
	
	A complementary comparison comes from degree-one cohomology. Palaisti \cite[Theorem~3.4 and Proposition~3.5]{PalaistiMetacyclic2026} shows that the subalgebras generated by degree-one mod-$2$ cohomology are isomorphic throughout the nonmodular metacyclic families considered here, and realizes them by inflation from a common quotient of order $16$. Applied to the groups in Proposition~\ref{prop:sharpness}, this uniformity shows that degree-one cup products do not resolve the ambiguity in the group order left by the first-layer transfer diagram. The higher class-field filtration recovers that missing information.
	
	The passage to a cyclotomic base changes the arithmetic problem. Mizusawa \cite{Mizusawa2005} studied finite nonabelian tower groups over real quadratic $\Z_2$-extensions with abelianization $C_2\times C_2$; his later metabelian examples \cite{Mizusawa2022} have abelianization $C_2^3$. These works relate the size of the nonabelian group to the finiteness of abelian Iwasawa modules. Our hypothesis $C_2\times C_{2^n}$, $n\ge2$, specifies a different case, in which one quadratic overfield is intrinsically distinguished. Recognition of metacyclicity through this overfield appears in Azizi--Rezzougui--Zekhnini \cite{AziziRezzouguiZekhnini2021}; transfer criteria for other rank-two abelianizations are developed in \cite{AziziRezzouguiZekhnini2022}. At finite level, Essahel--Dakkak--Mouhib \cite{EssahelDakkakMouhib2019} study which real quadratic fields have metacyclic Hilbert $2$-class towers. These are recognition and realization questions. The present paper instead asks what can be recovered after metacyclicity has been established.
	
	For this recovery problem to produce numerical arithmetic conclusions, one needs actual overfield modules. The rank formulas of Azizi--Mouhib \cite{AziziMouhib2001} and the unit-sensitive class-number relations of Kuroda \cite{Kuroda1950} explain the role of quadratic and biquadratic fields in obtaining those inputs. Chems-Eddin's study of real biquadratic Iwasawa modules \cite{ChemsJNT2026}, and its continuation with El Mamry \cite{ChemsElMamryII2026}, develop the associated rank and stabilization questions. We use the specific stable modules and metacyclicity theorem of Chems-Eddin--El Mamry \cite{ChemsElMamry2026} for the equal-quartic family in Section~\ref{sec:arithmetic}. That section states every imported conclusion used in the application, and explains the class-number relation that supplies the distinguished target order. The broader classification of real quadratic cyclotomic extensions with metacyclic tower group by Mizusawa--Mouhib \cite{MizusawaMouhib2027} places these examples in a general setting; our theorem does not require a new classification of the base fields.
	
	The resulting contribution is a recovery statement with a precise limit. The known first-layer calculations give a canonical initial quotient. Computing the derived subgroups at every higher degree shows that three apparently different abelian class fields are actually the same field, and that these common fields exhaust the full tower. Their orders recover its degree without selecting a chain or computing capitulation. Proposition~\ref{prop:sharpness} constructs an entire family with identical first-layer data and every permitted degree, showing why the additional layers are necessary. The arithmetic application then identifies an explicit triquadratic module whose order performs the first test.
	
	\section{Arithmetic background and the first-layer quotient}\label{sec:background}
	
	\subsection{Class groups and the two towers}
	For a number field $k$, the ordinary ideal class group is the group of fractional ideals modulo principal ideals, and $\Cl_2(k)$ is its Sylow $2$-subgroup. For a finite abelian $2$-group $A$, its $2$-rank is $\dim_{\mathbb F_2}(A/2A)$. Rank two specifies the number of generators, but does not determine the invariant factors or the order. The distinction is important because our degree formulas use module orders, not ranks alone.
	
	The ordinary Hilbert $2$-class field $\cH(k)$ is the largest abelian extension of $2$-power degree unramified at every finite and real place, and Artin reciprocity gives $\Gal(\cH(k)/k)\cong\Cl_2(k)$; see \cite{NSW}. The narrow class group uses only principal ideals with totally positive generators, and allows ramification at real places in the corresponding class field. The two conventions can differ at $2$. For example, $\Q(\sqrt3)$ has ordinary class number one, whereas $\Q(\sqrt3,i)/\Q(\sqrt3)$ is unramified at finite primes and becomes complex. The discriminants are $12$ and $144=12^2$, so its relative finite discriminant is trivial. Our real-place convention excludes this extension.
	
	Starting with $k^{[0]}=k$, define $k^{[r+1]}=\cH(k^{[r]})$. The union of these fields is $\cL(k)$. If $P=G(k)$ and $P^{(0)}=P$, $P^{(r+1)}=(P^{(r)})'$, then
	\[ k^{[r]}=\cL(k)^{P^{(r)}},\qquad \Gal(k^{[2]}/k)\cong P/P''. \]
	Thus the second Hilbert $2$-class field records a metabelian quotient, and a class group of an overfield records an abelianization of a subgroup. This is why overfield class groups can detect nonabelian structure in the original tower.
	
	The cyclotomic $\Z_2$-extension is a different tower. Let $\Q_\oo$ be the union of the maximal real subfields of the fields of $2$-power roots of unity, and put $k_\oo=k\Q_\oo$. Write $k_\nu$ for its unique subfield of degree $2^\nu$ over $k$. Iwasawa's formula \cite{Iwasawa1959,Washington} gives integers $\lambda,\mu\ge0$ and $\nu_0\in\Z$ such that, for all sufficiently large $\nu$,
	\[ |\Cl_2(k_\nu)|=2^{\lambda\nu+\mu2^\nu+\nu_0}. \]
	Greenberg's conjecture \cite{Greenberg1976} predicts $\lambda=\mu=0$ for totally real $k$, equivalently finiteness of $X(k_\oo)$. If $k/\Q$ is abelian, $\mu=0$ is already known by Ferrero--Washington \cite{FerreroWashington1979}. This applies to the quadratic and biquadratic fields in Section~\ref{sec:arithmetic}; their stable class-group calculations address the remaining boundedness question. The finiteness of the ordinary unramified module is an arithmetic input; it does not by itself give the size of the full nonabelian group $G(k_\oo)$.
	
	We will use stable finite-level class groups to describe such modules. The following standard result states exactly what is needed to justify that passage.
	
	\begin{proposition}\label{prop:stability} Let $k_\oo/k$ be a $\Z_2$-extension. Choose $n_0$ so that every prime ramified in $k_\oo/k$ is totally ramified in $k_\oo/k_{n_0}$. If $n\ge n_0$ and $|\Cl_2(k_{n+1})|=|\Cl_2(k_n)|$, then the norm maps identify $\Cl_2(k_\ell)$ with $\Cl_2(k_n)$ for every $\ell\ge n$, and $X(k_\oo)\cong\Cl_2(k_n)$. \end{proposition}
	\begin{proof} This is Fukuda's stabilization theorem \cite[Theorem~1(1)]{Fukuda1994}. Its mechanism can be seen after taking $k_n$ as base. Write $\Gamma=\langle\sigma\rangle$ and $\Lambda=\Z_2[[\sigma-1]]$. The class-field norm description used in that theorem gives a finitely generated $\Lambda$-submodule $Y\subseteq X=X(k_\oo)$ with $\Cl_2(k_{n+\ell})\cong X/\omega_\ell Y$, where $\omega_\ell=1+\sigma+\cdots+\sigma^{2^\ell-1}$ and $\omega_0=1$. Equal orders at the first two levels imply $(1+\sigma)Y=Y$. Since $1+\sigma$ lies in the maximal ideal $(2,\sigma-1)$ of $\Lambda$, Nakayama's lemma gives $Y=0$. The displayed description then identifies every subsequent class group with $X$. \end{proof}
	
	Equal ranks at two levels are a different assertion and do not replace the equality of orders in Proposition~\ref{prop:stability}. In the application, the cited arithmetic work supplies the required stabilization, including its ramification hypotheses.
	
	\subsection{Overfields, norms and capitulation}
	
	All derived subgroups of profinite groups are closed. For algebraic fields, unramifiedness is understood at the induced places, with the real-place convention above. The notation $X(E)$ always refers to the field $E$ itself; thus $X(k)$ and $X(k_\oo)$ are distinguished. For a finite unramified extension $E/B$, the two natural class-group maps have different roles: extension sends an ideal class of $B$ to its class in $E$, while norm sends classes in the opposite direction. Capitulation is the kernel of the extension map. The following dictionary identifies the group-theoretic counterpart of each.
	
	\begin{lemma}\label{lem:dictionary} Let $B$ be a number field or a totally real algebraic extension of $\Q$, and let $E/B$ be finite inside $\cL(B)$. Then $\cL(E)=\cL(B)$. If $H=\Gal(\cL(B)/E)$, then \[ X(E)\cong H/H',\qquad \cH(E)=\cL(B)^{H'}. \] The transfer $G(B)^{\ab}\to H^{\ab}$ corresponds to extension of ordinary ideal classes, or to its inverse-limit counterpart over a cyclotomic extension. The norm map corresponds to $H/H'\to G(B)/G(B)'$ induced by inclusion; if $E\subseteq\cH(B)$, its image is $H/G(B)'$. \end{lemma}
	\begin{proof} The extension $\cL(B)/E$ is unramified and pro-$2$, so it is contained in $\cL(E)$. Conversely, a finite Galois unramified $2$-extension of $E$ has normal closure over $B$ with a $2$-group as Galois group: take first the normal closure of $E/B$, whose group is a $2$-group, and then the compositum of the finitely many conjugate extensions. Unramifiedness is preserved by base change, composition and passage to subfields. The same operations preserve the condition at real places. The normal closure is therefore contained in $\cL(B)$. Taking unions proves equality. The abelian statements follow by fixing $H'$. The compatibility of extension with transfer and of norm with subgroup inclusion is Artin reciprocity; see \cite{Miyake1989,GruenbergWeiss2000,NSW}. Over a cyclotomic extension it follows by passing through finite levels and taking norm-compatible limits as in \eqref{eq:iwasawa}. \end{proof}
	
	The direction of these maps is worth making explicit. For $E\subseteq\cH(B)$, $[E:B]=[G:H]$ is read from the norm image $H/G'$. The capitulation map instead goes from $G/G'$ to the usually larger group $H/H'$. Knowing the source abelianization alone supplies no formula for the latter order.
	
	The transfer formulas needed in this paper have a short direct derivation.
	
	\begin{lemma}\label{lem:transferformula} Let $H$ be a subgroup of index two in a finite group $P$, and choose $z\in P\setminus H$. With $x^z=z^{-1}xz$, transfer is given by \[ \Ver_{P,H}(xP')=\begin{cases} xx^zH',&x\in H,\\ x^2H',&x\notin H.\end{cases} \] \end{lemma}
	\begin{proof} Use the transversal $\{1,z\}$. An element of $H$ fixes both cosets, producing the two conjugate factors in the first line. An element outside $H$ interchanges the cosets, so its two factors multiply to its square modulo $H'$. The product is independent of the transversal and factors through $P/P'$; see also \cite{Artin1929,Miyake1989,Mayer2016}. \end{proof}
	
	\subsection{The canonical first-layer quotient}
	For $G^{\ab}\cong C_2\times C_{2^n}$ there are three quadratic overfields because $G^{\ab}/2G^{\ab}\cong C_2^2$ has three nonzero homomorphisms to $C_2$. Their kernels lift to the three maximal subgroups of $G$. The first-layer transfer diagram records the abelianizations of these subgroups and the class-extension maps into them. Two diagrams are identified using compatible isomorphisms of their common source and their targets. In particular, capitulation kernels must be compared up to the same source automorphism, not by the spelling of a chosen generator.
	
	\begin{definition} For a finitely generated pro-$p$ group $P$ whose maximal subgroups have finite abelianization, put \[ N_1(P)=\bigcap_{H\text{ maximal in }P}H',\qquad S_1(P)=P/N_1(P). \] We call $S_1(P)$ the \emph{first-layer quotient}. Its first-layer transfer diagram consists of $P^{\ab}$, the maximal-subgroup abelianizations, and the transfers between them. \end{definition}
	
	The following elementary quotient observation isolates the information retained by those natural maps. Transfer naturality is standard; compare \cite{Mayer2016}.
	
	\begin{proposition}\label{prop:universal} The subgroup $N_1(P)$ is characteristic and open, and \[ P''\subseteq N_1(P)\subseteq P'. \] The finite quotient $S_1(P)$ preserves the first-layer transfer diagram through its natural maps. Every quotient preserving $P^{\ab}$ and all maximal-subgroup abelianizations through those maps surjects naturally onto $S_1(P)$. Moreover, \[ N_1(S_1(P))=1,\qquad S_1(S_1(P))\cong S_1(P). \] \end{proposition}
	\begin{proof} There are finitely many maximal subgroups, all open of index $p$. Each $H'$ is open because $H/H'$ is finite, so their intersection is open. Automorphisms permute the $H'$, and $P'\subseteq H$ gives $P''\subseteq H'\subseteq P'$. \par Write $N=N_1(P)$. The maximal subgroups of $P/N$ are $H/N$, with derived subgroups $H'/N$. Thus the natural maps identify the source and targets of the first-layer diagram; transfer naturality identifies its arrows. If a normal subgroup $D$ defines another quotient preserving the source, then $D\subseteq P'$. Preservation of $H/H'$ then requires $D\subseteq H'$ for every maximal $H$, so $D\subseteq N$. Finally, \[ \bigcap_H(H/N)'=(\bigcap_H H')/N=1, \] which proves the last assertions. \end{proof}
	
	In arithmetic terms, $N_1(G)=\bigcap_H H'$ fixes the compositum of the three fields $\cH(E)$ above the quadratic overfields. The inclusions $G''\subseteq N_1(G)\subseteq G'$ become \[ \cH(B)\subseteq\cL(B)^{N_1(G)}\subseteq\cH(\cH(B)). \] Thus the first-layer quotient is an intermediate part of the second Hilbert $2$-class field. It need not contain every part of that field. Proposition~\ref{prop:firstquotient} identifies this compositum with the class field above the distinguished quadratic overfield in the metacyclic case.
	
	This universal property concerns quotients of a given group and their natural maps. Proposition~\ref{prop:sharpness} gives explicit failures of such determination, even within the metacyclic class.
	
	Finiteness of the abelianization alone would not suffice, even for abstract metacyclic pro-$2$ groups. For example, let $C_{2^n}$ act on $\Z_2$ by inversion. The infinite semidirect product $\Z_2\rtimes C_{2^n}$ has abelianization $C_2\times C_{2^n}$, but its distinguished maximal subgroup is the infinite abelian group $\Z_2\times C_{2^{n-1}}$. The finite overfield modules in the next lemma exclude exactly this kind of possibility.
	
	\begin{lemma}\label{lem:finite} Let $P$ be metacyclic pro-$2$. If $P^{\ab}$ and $H^{\ab}$ are finite for every maximal subgroup $H$, then $P$ is finite. \end{lemma}
	\begin{proof} Choose a closed normal procyclic subgroup $A$ with procyclic quotient $P/A$. This quotient is a quotient of $P^{\ab}$ and is finite. If $A$ is infinite, then $A\cong\Z_2$. Its conjugation action factors through $P/A$, and a finite $2$-subgroup of $\Z_2^\times$ has order at most two. The preimage of the action kernel has index at most two, centralizes $A$, and is abelian because its quotient by $A$ is cyclic. It is infinite, contradicting either finiteness of $P^{\ab}$ or that of a maximal-subgroup abelianization. Thus $A$ and $P$ are finite. \end{proof}
	
	\section{Metacyclic presentations and the first layer}\label{sec:metacyclic}
	
	A finite $2$-group is metacyclic if it has a cyclic normal subgroup with cyclic quotient. The modular group of order $2^\ell$, $\ell>3$, has presentation
	\[ \langle x,y\mid x^2=y^{2^{\ell-1}}=1,\ [x,y]=y^{2^{\ell-2}}\rangle. \]
	For the groups considered below, a presentation adapted to the cyclic normal subgroup has the form \[ P=\langle a,b\mid a^{2^\alpha}=1,\ b^{2^n}=a^c,\ b^{-1}ab=a^u\rangle. \] The odd integer $u$ describes the action on $\langle a\rangle$, while $c$ describes the power relation for a lift of the cyclic quotient. Consistency requires $u^{2^n}\equiv1\pmod{2^\alpha}$ and $(u-1)c\equiv0\pmod{2^\alpha}$. These congruences are also sufficient: start with $C_{2^\alpha}\rtimes\Z$, on which an infinite cyclic generator $b$ acts by $a\mapsto a^u$, and quotient by the central element $b^{2^n}a^{-c}$. Its nonzero image in $\Z$ ensures that the quotient retains the full subgroup $C_{2^\alpha}$ and has cyclic quotient of order $2^n$. Every element then has a unique normal form $b^j a^i$, $0\le j<2^n$, $0\le i<2^\alpha$. Multiplication uses $ab=ba^u$ and the power relation for $b$, so it reduces to arithmetic modulo $2^\alpha$. King's presentations \cite{King1973} and Xu--Zhang's classification \cite{XuZhang2006} treat this general problem and the identifications between presentations.
	
	Following \cite{AabounZekhnini2023}, metacyclic-nonmodular means nonabelian metacyclic with this modular family excluded. For abelianization $C_2\times C_{2^n}$, $n\ge2$, the following presentations suffice:
	\begin{align*} \mathrm{I}:\quad&a^{2^\alpha}=b^{2^n}=1, &b^{-1}ab&=a^{-1},&\alpha&\ge2;\\ \mathrm{II}:\quad&a^{2^\alpha}=1,\ b^{2^n}=a^{2^{\alpha-1}}, &b^{-1}ab&=a^{-1},&\alpha&\ge3;\\ \mathrm{III}:\quad&a^{2^\alpha}=b^{2^n}=1, &b^{-1}ab&=a^{-1+2^s},&2&\le s\le\alpha-1;\\ \mathrm{IV}:\quad&a^{2^\alpha}=1,\ b^{2^n}=a^{2^{\alpha-1}}, &b^{-1}ab&=a^{-1+2^s},&2&\le s\le\alpha-2. \end{align*}
	The displayed ranges are supplemented by $\alpha\le n+s$ in families~III--IV. This is precisely the further condition needed for the indicated order: $u^{2^n}\equiv1\pmod{2^\alpha}$ is equivalent to that inequality by \eqref{eq:LTE}, and the condition on $c$ is automatic. The construction above therefore gives $a$ of order $2^\alpha$ and cyclic quotient of order $2^n$, so $|P|=2^{n+\alpha}$. The labels refer to presentation families, and are not used as a disjoint classification of isomorphism classes. Some presentations in different rows are isomorphic; Lemma~\ref{lem:endpoint} makes one such overlap explicit. Boundary family~III means $s=\alpha-1$. The normalization of the action exponent can be justified directly. Suppose $u=-1+\kappa2^s$, where $\kappa$ is odd and $2\le s\le\alpha-1$. The units congruent to $1$ modulo $2^s$, considered modulo $2^\alpha$, form a group of order $2^{\alpha-s}$. The same valuation calculation as in \eqref{eq:LTE} shows that both $-u=1-\kappa2^s$ and $1-2^s$ have this order, so each generates that group. Hence there is an odd integer $h$ such that $(-u)^h\equiv1-2^s\pmod{2^\alpha}$, and therefore $u^h\equiv-1+2^s\pmod{2^\alpha}$. Replacing $b$ by $b^h$ preserves its generating image in the cyclic quotient and gives the normalized action. It also preserves the power relation: if $b^{2^n}=a^c$ with $c\in\{0,2^{\alpha-1}\}$, then $(b^h)^{2^n}=a^{ch}=a^c$, since $h$ is odd.
	
	The action in families~I--II is exact inversion. In families~III--IV, $s=v_2(u+1)$ measures its deviation from inversion, where $v_2(z)$ is the exponent of $2$ dividing a nonzero integer $z$. This parameter becomes visible when one passes to $b^2$: the action is then multiplication by $u^2$, and $v_2(u^2-1)=s+1$. Each further squaring raises that valuation by one. This is the source of the common higher derived subgroups and, ultimately, of the plateau in Theorem~\ref{thm:main}.
	
	Put $u=-1$ in families~I--II and $u=-1+2^s$ in families~III--IV. In all cases $P'=\langle a^2\rangle$, and the maximal subgroups are
	\begin{equation} H_1=\langle b,a^2\rangle,\qquad H_2=\langle ab,a^2\rangle,\qquad H_3=\langle a,b^2\rangle. \end{equation}
	The image of $H_3$ is the unique noncyclic index-two subgroup of $P^{\ab}$. Hence $H_3$ is characteristic.
	
	The next formulas are those of \cite[Corollary~3.3]{AabounZekhnini2023}; the proof is recalled to fix the generators used in the higher-layer calculation.
	
	\begin{lemma}\label{lem:firstcomm} For the presentations above, \[ H_1'=H_2'=\langle a^4\rangle,\qquad H_3'=\begin{cases} 1,&\mathrm{I},\mathrm{II},\\ \langle a^{2^{s+1}}\rangle,&\mathrm{III},\mathrm{IV}. \end{cases} \] In particular, $N_1(P)=H_3'$. \end{lemma}
	\begin{proof} Since $v_2(u-1)=1$, $P'=\langle a^{u-1}\rangle=\langle a^2\rangle$. The commutators of $a^2$ with $b$ and $ab$ generate $\langle a^{2(u-1)}\rangle=\langle a^4\rangle$. The group $H_3'$ is generated by $a^{u^2-1}$. It is trivial for $u=-1$; otherwise \[ u^2-1=2^{s+1}(-1+2^{s-1}), \] whose final factor is odd. Thus $H_3'\subseteq H_1'=H_2'$, giving the intersection formula. \end{proof}
	
	\begin{proposition}\label{prop:firstquotient} For a group $P$ as above, \[ S_1(P)=P/H_3',\qquad |S_1(P)|=2|H_3^{\ab}|. \] It equals $P$ precisely when $H_3$ is abelian, namely in families~I, II and boundary~III. Otherwise it has the boundary-III presentation with parameters $n$, $\alpha_0=s+1$. Thus an actual first-layer quotient in the metacyclic-nonmodular class has a presentation in one of these three boundary families. \end{proposition}
	\begin{proof} The order formula and the condition for equality follow from Lemma~\ref{lem:firstcomm}. In the remaining cases quotient by $\langle a^{2^{s+1}}\rangle$. The image of $a$ has order $2^{s+1}$ and conjugation exponent $-1+2^s$. In family~IV, $\alpha-1\ge s+1$, so the power relation becomes $b^{2^n}=1$. This is the asserted boundary presentation. The last statement also follows from the idempotence in Proposition~\ref{prop:universal}. \end{proof}
	
	The same conclusion holds when only $S_1(P)$, rather than $P$, is assumed metacyclic-nonmodular. Indeed, idempotence makes its distinguished maximal subgroup abelian, so $H_3'\subseteq N_1(P)$. The reverse inclusion is automatic. This gives $N_1(P)=H_3'$ and the same quotient order.
	
	\begin{lemma}\label{lem:endpoint} In families~III--IV, \begin{equation}\label{eq:bound} \alpha\le n+s. \end{equation} If $\alpha=n+s$, the presentations in families~III and~IV define isomorphic groups. \end{lemma}
	\begin{proof} Conjugation by $b^{2^n}$ fixes $a$, so $u^{2^n}\equiv1\pmod{2^\alpha}$. The elementary $2$-adic valuation identity gives \begin{equation}\label{eq:LTE} v_2(u^{2^j}-1)=s+j\qquad(j\ge1). \end{equation} For example, the case $j=1$ follows from the factorization in Lemma~\ref{lem:firstcomm}; each subsequent squaring adds one because $u^{2^j}\equiv1\pmod8$. Taking $j=n$ proves \eqref{eq:bound}. \par Now suppose $\alpha=n+s$ and begin with family~IV. Set $c=ba$. Then $c^{-1}ac=a^u$ and \[ c^{2^n}=b^{2^n}a^S,\qquad S=1+u+\cdots+u^{2^n-1}. \] As $v_2(u-1)=1$, \eqref{eq:LTE} gives $v_2(S)=n+s-1=\alpha-1$. Hence $a^S=a^{2^{\alpha-1}}$ and $c^{2^n}=1$. The generators $a,c$ give the family-III presentation with the same parameters. Both groups have order $2^{n+\alpha}$, proving the isomorphism. \end{proof}
	
	\section{Common class fields and recovery of tower degrees}\label{sec:higher}
	
	For a finite group $P$ with $P^{\ab}\cong C_2\times C_{2^n}$, define
	\begin{equation} \mathcal U_j(P)=\{U:P'\subseteq U\subseteq P,\ [P:U]=2^j\},\quad N_j(P)=\bigcap_{U\in\mathcal U_j(P)}U',\quad S_j(P)=P/N_j(P). \end{equation}
	These definitions are independent of a presentation. The case $j=1$ is the first-layer quotient. All groups in $\mathcal U_j(P)$ are normal. The argument of Proposition~\ref{prop:universal}, with $\mathcal U_j(P)$ in place of the maximal subgroups, shows that $S_j(P)$ is the smallest quotient preserving $P^{\ab}$ and every degree-$2^j$ transfer target and map through the natural maps. It also gives $N_j(S_j(P))=1$. Thus the filtration below consists of universal quotients at the successive degrees.
	
	Lemma~\ref{lem:dictionary} translates the equality of class fields in Theorem~\ref{thm:main} into a subgroup calculation. Two overfields $E_1,E_2$ have the same maximal unramified abelian extension precisely when their corresponding subgroups have the same derived subgroup. We therefore list the three subgroups at a given degree and calculate these derived subgroups, rather than infer equality of fields from equality of class numbers.
	
	\begin{lemma}\label{lem:highercomm} Let $P$ have one of the presentations in Section~\ref{sec:metacyclic}. For $2\le j\le n$, $\mathcal U_j(P)$ consists of \begin{equation}\label{eq:threeU} \begin{split} U_{1,j}&=\langle a,b^{2^j}\rangle,\\ U_{2,j}&=\langle ab^{2^{j-1}},a^2\rangle,\\ U_{3,j}&=\langle b^{2^{j-1}},a^2\rangle. \end{split} \end{equation} Their derived subgroups coincide. They are trivial in families~I--II, and in families~III--IV they are \begin{equation} U_{1,j}'=U_{2,j}'=U_{3,j}'=\langle a^{2^{s+j}}\rangle. \end{equation} \end{lemma}
	\begin{proof} Write $A=P/P'=\langle\bar a\rangle\times\langle\bar b\rangle$. An index-$2^j$ subgroup of $A$ either contains $\bar a$, in which case it is $\langle\bar a,\bar b^{2^j}\rangle$, or meets $\langle\bar a\rangle$ trivially. In the latter case it projects isomorphically onto $\langle\bar b^{2^{j-1}}\rangle$ and is the graph of one of its two homomorphisms to $\langle\bar a\rangle$. These are the other two subgroups displayed in \eqref{eq:threeU}. Taking full preimages proves the list. \par The first derived group is generated by $a^{u^{2^j}-1}$. The third is generated by $a^{2(u^{2^{j-1}}-1)}$, and the second has the same generator up to a unit power: $ab^{2^{j-1}}$ acts on $\langle a^2\rangle$ as $b^{2^{j-1}}$ does. These elements lie in the normal cyclic subgroup $\langle a\rangle$, and their cyclic subgroups are invariant under the indicated generators, so they generate the full derived groups. For $u=-1$ they vanish. Otherwise \eqref{eq:LTE}, with $j-1\ge1$, gives valuation $s+j$ in all three cases. \end{proof}
	
	For $j=2$, this is exactly the index-four calculation of \cite[Corollary~3.3]{AabounZekhnini2023}. The point of the uniform expression is that it gives both equality of the fields fixed by these subgroups and the whole filtration, not just a table of quartic target orders.
	
	\begin{proposition}\label{prop:filtration} Let $P$ be as in Section~\ref{sec:metacyclic}, and write $|H_3^{\ab}|=2^{n+t}$. Then \[ t=\alpha-1\quad\hbox{in families I--II and boundary III}, \qquad t=s\quad\hbox{in the remaining cases}. \] Set $k=\alpha-t$. Then $1\le k\le n$, and, for $1\le j\le n$, \begin{equation}\label{eq:filtrationkernel} N_j(P)=\langle a^{2^{t+j}}\rangle,\qquad |S_j(P)|=2^{n+t+\min\{j,k\}}. \end{equation} The kernels are characteristic and form a decreasing chain, with $N_n(P)=1$. For $j\ge2$, $N_j(P)=U'$ for every $U\in\mathcal U_j(P)$, and \begin{equation}\label{eq:allorders} |U^{\ab}|=2^{n+t-\max\{j-k,0\}}. \end{equation} \end{proposition}
	\begin{proof} In the first three cases $H_3$ is abelian of index two, so $t=\alpha-1$ and $k=1$. All the subgroups at degrees $2^j$, $j\ge2$, are then abelian by Lemma~\ref{lem:highercomm}; the expression for $N_j(P)$ is also trivial since $t+j\ge\alpha$. \par Otherwise $|H_3'|=2^{\alpha-s-1}$, and \[ |H_3^{\ab}|=2^{n+\alpha-1}/2^{\alpha-s-1}=2^{n+s}. \] Thus $t=s$ and $k=\alpha-s$. Lemmas~\ref{lem:firstcomm} and~\ref{lem:highercomm} give the kernels for every $j$. Their orders are $2^{\max\{k-j,0\}}$, yielding \eqref{eq:filtrationkernel} and \eqref{eq:allorders}. The bound $k\le n$ is \eqref{eq:bound}. Automorphisms preserve the sets $\mathcal U_j(P)$, so their derived-subgroup intersections are characteristic. Finally $t+n\ge\alpha$ gives $N_n(P)=1$. \end{proof}
	
	The parameters now have a direct interpretation. The source module gives $n$, and the distinguished quadratic module gives $t$. The full cyclic-normal parameter $\alpha$ is hidden at the first layer. The integer $k=\alpha-t$ measures how many successive factors of two the canonical quotient can gain before it becomes the full group. Formula~\eqref{eq:allorders} turns this missing parameter into a visible length in a sequence of abelian module orders.
	
	\begin{proof}[Proof of Theorem~\ref{thm:main}] Under the Galois correspondence, the fields $E$ of degree $2^j$ in $\cH(B)$ correspond to $\mathcal U_j(G)$. By Lemma~\ref{lem:dictionary}, $\cH(E)=\cL(B)^{U'}$. For $j\ge2$ all three derived subgroups equal $N_j(G)$, so the three fields coincide with $\cL(B)^{N_j(G)}$. At $j=1$, Proposition~\ref{prop:firstquotient} identifies the compositum of the three class fields with $\cH(K_B)=\cL(B)^{N_1(G)}$. Proposition~\ref{prop:filtration} proves the Galois and filtration statements and \eqref{eq:mainfields}--\eqref{eq:mainorders}. In fact each $\cM_j$ is invariant under every automorphism of $\cL(B)/B$. \par The integer $k$ is the last index for which the target order is $2^{n+t}$. If $k=1$, the quartic order is $2^{n+t-1}$; if $k\ge2$, it is $2^{n+t}$. Moreover $N_1(G)=1$ exactly when $k=1$. This proves \eqref{eq:quarticintro}. Taking $j=n$ in \eqref{eq:mainorders} gives $w_n=2^{t+k}$, and \eqref{eq:mainrecovery} follows. \end{proof}
	
	For example, if $n=3$ and $t=2$, the distinguished quadratic module has order $32$. The three possible order sequences are displayed in Table~\ref{tab:plateau}. Here $w_2$ may be read from any quartic field and $w_3$ from any degree-eight field in $\cH(B)$; the choices need not be nested.
	
	\begin{table}[ht]
		\centering
		\caption{Reading the full degree from the plateau when $n=3$, $t=2$.}\label{tab:plateau}
		\begin{tabular}{@{}ccc@{}}
			\toprule
			$k$&$(w_1,w_2,w_3)$&$[\cL(B):B]$\\
			\midrule
			$1$&$(32,16,8)$&$64$\\
			$2$&$(32,32,16)$&$128$\\
			$3$&$(32,32,32)$&$256$\\
			\bottomrule
		\end{tabular}
	\end{table}
	
	When $B=k_\oo$, this canonicity also respects the cyclotomic Galois action. Indeed, every $k$-automorphism of an algebraic closure preserves $B$ and its maximal unramified pro-$2$ extension, so $\cL(B)/k$ is Galois. Conjugation preserves the characteristic subgroups $N_j(G)$; consequently every $\cM_j/k$ is Galois. The filtration is therefore compatible with the arithmetic action of $\Gal(k_\oo/k)$, although its degree formulas require only the underlying group and module orders.
	
	Thus $\cM_j$ is the compositum of the maximal unramified abelian extensions of the degree-$2^j$ fields in the first abelian layer. For $j\ge2$, taking the compositum is unnecessary: any one of those fields gives the same result. One should distinguish these degrees over $B$ from the degrees of successive Hilbert class fields. The filtration is indexed by intermediate fields inside $\cH(B)$.
	
	\begin{corollary}\label{cor:endpoint} Assume Theorem~\ref{thm:main}, with $t\ge2$. If $k=n$, equivalently $|G|=2^{2n+t}$, then the isomorphism class of $G$ is determined by $n,t$: \[ G\cong\langle a,b\mid a^{2^{n+t}}=b^{2^n}=1,\ b^{-1}ab=a^{-1+2^t}\rangle. \] In particular, if $n=2$ and a quartic target has order $|X(K_B)|$, this presentation gives the full tower group. \end{corollary}
	\begin{proof} As $n\ge2$, $k=n$ excludes the cases with abelian $H_3$. The remaining presentations have $s=t$, $\alpha=n+t$. Lemma~\ref{lem:endpoint} identifies the two possible power relations. For $n=2$, the quartic test distinguishes $k=1$ from $k=2$. \end{proof}
	
	\begin{proposition}\label{prop:sharpness} Fix $n\ge2$ and $t\ge2$. For each $1\le r\le n$, let \[ P_r=\langle a,b\mid a^{2^{t+r}}=b^{2^n}=1,\ b^{-1}ab=a^{-1+2^t}\rangle. \] Then: \begin{enumerate} \item $|P_r|=2^{n+t+r}$, and these $n$ pairwise nonisomorphic groups have the same complete first-layer transfer diagram under natural identifications. Their common first-layer quotient is $P_1$. \item Their higher canonical quotients satisfy $S_j(P_r)\cong P_{\min\{j,r\}}$ for $1\le j\le n$. \item If $n\ge3$, there are also two nonisomorphic groups of order $2^{n+t+2}$ with the same complete first-layer transfer diagram, obtained by taking $\alpha=t+2$, $s=t$ in families~III and~IV. \end{enumerate} \end{proposition}
	\begin{proof} The automorphism $a\mapsto a^{-1+2^t}$ has order dividing $2^n$ on $C_{2^{t+r}}$, because \eqref{eq:LTE} gives $n+t\ge t+r$. Thus $P_r$ is the indicated semidirect product of order $2^{n+t+r}$, and $P_r'=\langle a^2\rangle$ gives $P_r^{\ab}\cong C_2\times C_{2^n}$. Its first-layer kernel is $\langle a^{2^{t+1}}\rangle$, and the quotient is $P_1$, compatibly with the generator images. Proposition~\ref{prop:universal} identifies all transfer targets, maps and kernels. Distinct orders prove nonisomorphism. The same calculation with $N_j(P_r)$ proves (2). \par For (3), both presentations are consistent extensions of $C_{2^\alpha}$ by $C_{2^n}$; in family~IV, the central power $a^{2^{\alpha-1}}$ is fixed by the action. Both have first-layer quotient $P_1$. To distinguish them, consider elements whose image in the abelianization has order $2^n$. They are precisely $b^j a^i$ with $j$ odd. For such a $j$, $u^j\equiv-1\pmod4$ and \[ v_2\left(1+u^j+\cdots+(u^j)^{2^n-1}\right)=n+t-1\ge t+2=\alpha. \] Indeed $v_2(u^j+1)=t$ and $v_2(u^j-1)=1$, so the assertion follows as in \eqref{eq:LTE}. Hence the $2^n$th power of $b^j a^i$ is $1$ in family~III and $a^{2^{\alpha-1}}$ in family~IV. Such elements have order $2^n$ in the former and $2^{n+1}$ in the latter. This is an isomorphism invariant, proving nonisomorphism. \end{proof}
	
	The restriction $n\ge3$ in part~(3) is essential: for $n=2$, $\alpha=t+2$ is the endpoint of Lemma~\ref{lem:endpoint}, and those two presentations are isomorphic. Part~(1), which also covers $n=2$, proves that all the order bounds in Theorem~\ref{thm:main} are sharp among abstract groups with a common first-layer diagram. It does not assert that every such group is realized by the quadratic arithmetic family below.
	
	\begin{example}\label{ex:small} Take $n=t=2$. The two groups \[ P_1=\langle a,b\mid a^8=b^4=1,\ b^{-1}ab=a^3\rangle, \qquad P_2=\langle a,b\mid a^{16}=b^4=1,\ b^{-1}ab=a^3\rangle \] have orders $32$ and $64$ and the same first-layer transfer diagram. The distinguished quadratic target has order $16$ in both. Every quartic target has order $8$ for $P_1$ and order $16$ for $P_2$. Thus the quartic layer detects an extra factor of two in the full degree even when all quadratic transfer data agree. \end{example}
	
	\section{Application to a family of real quadratic fields}\label{sec:arithmetic}
	
	Let $\eta\in\{1,2\}$, and let $q,r,s$ be distinct primes satisfying
	\[ q\equiv3\pmod4,\qquad r\equiv s\equiv5\pmod8. \]
	Put
	\begin{equation} \begin{aligned} F&=\Q(\sqrt{\eta qrs}),&\qquad K&=\Q(\sqrt{\eta q},\sqrt{rs}),\\ K'&=\Q(\sqrt r,\sqrt{\eta qs}),&\qquad K''&=\Q(\sqrt s,\sqrt{\eta qr}). \end{aligned} \end{equation}
	For positive squarefree $d$, let $h_2(d)$ denote the $2$-part of the ordinary class number of $\Q(\sqrt d)$, and let $\epsilon_d>1$ be its fundamental unit. Write $\delta_{a,b}$ for the Kronecker delta. We impose the hypotheses of \cite[Theorem~3.10]{ChemsElMamry2026}:
	\begin{gather} \left(\frac qr\right)=\left(\frac qs\right)=(-1)^{\delta_{\eta,2}},\qquad \left(\frac rs\right)=1,\qquad N_{\Q(\sqrt{rs})/\Q}(\epsilon_{rs})=1,\label{eq:arithhyp}\\ 2^{\delta_{\eta,1}}q(\gamma-1)\notin\{z^2:z\in\Z\},\qquad \epsilon_{\eta qrs}=\gamma+\gamma'\sqrt{\eta qrs}.  \end{gather}
	In addition, assume
	\begin{equation}\label{eq:equal} \left(\frac rs\right)_4=\left(\frac sr\right)_4=1,\qquad h_2(\eta qrs)=2^{m+1},\quad h_2(rs)=2^t,\quad m,t\ge2. \end{equation}
	Here the rational quartic symbol, when the Legendre symbol is $1$, is $1$ or $-1$ according as the numerator is a fourth power or a square but not a fourth power modulo the denominator.
	
	The next proposition collects the imported arithmetic results in the notation used here. Its conclusions, rather than any additional assertion about the full group order, are the inputs to the recovery theorem.
	
	\begin{proposition}\label{prop:arithinput} Assume \eqref{eq:arithhyp}--\eqref{eq:equal}. The extensions $K_\oo$, $K'_\oo$, $K''_\oo$ are the three distinct unramified quadratic extensions of $F_\oo$. The subgroup corresponding to $K_\oo$ has noncyclic image in $X(F_\oo)$; equivalently, every quartic subfield of $\cH(F_\oo)/F_\oo$ contains $K_\oo$. Moreover, \[ \begin{aligned} X(F_\oo)&\cong C_2\times C_{2^m},& |X(K_\oo)|&=2^{m+t},\\ X(K'_\oo)&\cong C_2\times C_{2^m},&X(K''_\oo)&\cong C_2\times C_{2^m}. \end{aligned} \] Finally, $G(F_\oo)$ is a metacyclic pro-$2$ group. \end{proposition}
	\begin{proof} The overfield identification, source module and metacyclicity are \cite[Theorem~3.10 and its proof]{ChemsElMamry2026}. The order of $X(K_\oo)$ is Corollary~3.9 there. Corollary~3.12 gives $X(K'_\oo)$, and symmetry in $r,s$ gives the statement for $K''_\oo$. These results are obtained from finite-level class groups, with stabilization justified by the ramification and unit calculations in that work and Proposition~\ref{prop:stability}. \end{proof}

	The order formula for the distinguished module has a useful arithmetic explanation. Write $U_M=\mathcal O_M^\times$ for a unit group and put $k_1=\Q(\sqrt{\eta q})$, $k_2=\Q(\sqrt{rs})$, $k_3=F$. For the real biquadratic field $K$, Kuroda's class-number formula \cite{Kuroda1950} reads \[ h(K)=\frac{Q(K)}4 h(k_1)h(k_2)h(k_3),\qquad Q(K)=[U_K:U_{k_1}U_{k_2}U_{k_3}]. \] The factor $Q(K)$ records units of $K$ not generated by units of its quadratic subfields. Under our hypotheses, the unit calculation in \cite[Corollary~3.9]{ChemsElMamry2026} gives $Q(K)=2$ and $h_2(k_1)=1$. Taking $2$-parts therefore gives \[ h_2(K)=\tfrac12 h_2(rs)h_2(\eta qrs)=2^{m+t}. \] The same order occurs at the next cyclotomic level, with the required ramification, so Proposition~\ref{prop:stability} identifies it with $|X(K_\oo)|$. The class-group rank calculations behind the other inputs are related to \cite{AziziMouhib2001,ChemsJNT2026}; rank stability alone would not yield this order formula.
	
	Proposition~\ref{prop:arithinput} supplies finite abelianizations for all three maximal subgroups of $G(F_\oo)$. Together with metacyclicity, Lemma~\ref{lem:finite} proves that this full group is finite. It cannot be abelian: that would give order $2^{m+1}$ and maximal-subgroup order $2^m$, smaller than $|X(K_\oo)|=2^{m+t}$. It cannot be modular either: a modular group with this abelianization has order $2^{m+2}$, so any maximal-subgroup abelianization has order at most $2^{m+1}$, again too small because $t\ge2$. Its distinguished subgroup is therefore the $H_3$ of Section~\ref{sec:metacyclic}, and all hypotheses of the recovery theorem have been checked.
	
	\begin{theorem}\label{thm:arithmetic} Under \eqref{eq:arithhyp}--\eqref{eq:equal}, put $B=F_\oo$, $G=G(B)$ and $G_1=\Gal(\cH(K_\oo)/B)$. Then \begin{equation} |G_1|=2^{m+t+1},\qquad 2^{m+t+1}\le[\cL(B):B]\le2^{2m+t}. \end{equation} The quotient $G_1$ has a presentation in family~I, II or boundary~III, with $n=m$ and $\alpha_0=t+1$. \par For every $2\le j\le m$, the three degree-$2^j$ fields $E/B$ inside $\cH(B)$ have the same field $\cH(E)$. For a single integer $1\le k\le m$, their degrees and modules satisfy \[ [\cH(E):B]=2^{m+t+\min\{j,k\}},\qquad |X(E)|=2^{m+t-\max\{j-k,0\}},\qquad |G|=2^{m+t+k}. \] In particular, the quartic test \eqref{eq:quarticintro} and recovery formula \eqref{eq:mainrecovery} apply with $K_B=K_\oo$. If $k=m$, the full group has the unique presentation of Corollary~\ref{cor:endpoint} with $n=m$. \end{theorem}
	\begin{proof} Proposition~\ref{prop:arithinput} identifies the distinguished maximal subgroup with $\Gal(\cL(B)/K_\oo)$. Theorem~\ref{thm:main} gives every degree formula and bound. Proposition~\ref{prop:firstquotient} gives the boundary presentations. Their distinguished maximal subgroups are abelian of order $2^{m+\alpha_0-1}$; comparison with $|X(K_\oo)|=2^{m+t}$ gives $\alpha_0=t+1$. The endpoint assertion is Corollary~\ref{cor:endpoint}. \end{proof}
	
	In the full group, $k=1$ means that $G=G_1$. Otherwise $s=t$ and $\alpha=t+k$, with $t+2\le\alpha\le t+m$. These are restrictions on admissible presentations; the overlap described in Lemma~\ref{lem:endpoint} is retained. Proposition~\ref{prop:sharpness} shows that no smaller general range of orders follows solely from first-layer transfer data in the abstract metacyclic class.
	
	\begin{corollary}\label{cor:typeII} Assume Theorem~\ref{thm:arithmetic} and $m\ge t+2$. Then $G_1$ has a family-II presentation if and only if \begin{equation} \expn X(K_\oo)=2^m. \end{equation} In this case $G=G_1$. The distinguished modules are \[ X(K_\oo)\cong \begin{cases} C_{2^t}\times C_{2^m},&\text{family II},\\ C_{2^{t+1}}\times C_{2^{m-1}},&\text{family I or boundary III}. \end{cases} \] \end{corollary}
	\begin{proof} Write $\alpha_0=t+1$. In family~I or boundary~III, the distinguished subgroup is $C_{2^{\alpha_0}}\times C_{2^{m-1}}$. In family~II it is abelian, generated by $a$ and $c=b^2$, with relations \[ a^{2^{\alpha_0}}=1,\qquad c^{2^{m-1}}=a^{2^{\alpha_0-1}}. \] Its invariant factors are $2^{\min\{\alpha_0-1,m-1\}}$ and $2^{\max\{\alpha_0,m\}}$, by the Smith normal form of the two relations. For $m\ge t+2$, these give the displayed alternatives, whose exponents differ. A proper first-layer quotient of $G$ has a boundary-III presentation by Proposition~\ref{prop:firstquotient}; it cannot be family~II, as their distinguished transfer kernels below differ. Thus the family-II alternative forces $G=G_1$. \end{proof}
	
	\begin{proposition}\label{prop:capitulation} In Theorem~\ref{thm:arithmetic}, choose $c=aG'$ and $d=bG'$ from a presentation compatible with $G_1$, so \[ X(B)=\langle c\rangle\times\langle d\rangle,\qquad |c|=2,\quad |d|=2^m. \] The kernel of $X(B)\to X(K_\oo)$ is \[ \begin{array}{c|c} \text{presentation of }G_1&\text{kernel}\\\hline \mathrm{I}&\langle c,d^{2^{m-1}}\rangle\\ \mathrm{II}&\langle c\rangle\\ \text{boundary III}&\langle d^{2^{m-1}}\rangle. \end{array} \] These three possibilities are distinct up to automorphisms of $X(B)$. \end{proposition}
	\begin{proof} Transfer to $H_3$ sends $aG'$ to $a^{1+u}H_3'$ and $bG'$ to $b^2H_3'$. Transfer is preserved in $G_1$, so compute there, where $H_3'=1$. In family~I both $a$ and $b^{2^{m-1}}$ transfer to $1$. In family~II, $a$ transfers to $1$ but $b^{2^{m-1}}$ transfers to $a^{2^{\alpha_0-1}}\ne1$. In boundary~III the reverse holds: $a$ transfers to $a^{2^{\alpha_0-1}}\ne1$ and $b^{2^{m-1}}$ transfers to $1$. The relations in the abelian target give the stated full kernels. These are also the distinguished transfer formulas of \cite[Theorem~4.5]{AabounZekhnini2023}. \par The first kernel has order four. Among the other two, $\langle d^{2^{m-1}}\rangle$ is the unique subgroup of order two in the characteristic subgroup $2X(B)$, while $\langle c\rangle$ is not contained in $2X(B)$. Hence no source automorphism identifies the three possibilities. \end{proof}
	
	\section{Explicit determination of tower groups}
	
	Table~\ref{tab:examples} uses five rows of the stable ordinary class-group data in \cite[Table~6]{ChemsElMamry2026}. The values of $t$ follow from $|X(K_\oo)|=2^{m+t}$. Parentheses abbreviate invariant factors. The last three columns are consequences of Theorem~\ref{thm:arithmetic}, not new computations of class groups.
	
	\begin{table}[ht]
		\centering\small
		\caption{Stable inputs and resulting tower-degree candidates.}\label{tab:examples}
		\begin{tabular}{@{}lccccc@{}}
			\toprule
			$(\eta,q,r,s)$&$X(F_\oo)$&$X(K_\oo)$&$(m,t)$&$|G_1|$&possible $|G|$\\
			\midrule
			$(1,131,53,13)$&$(2,4)$&$(2,8)$&$(2,2)$&$32$&$32,64$\\
			$(1,367,61,13)$&$(2,8)$&$(4,8)$&$(3,2)$&$64$&$64,128,256$\\
			$(1,491,61,13)$&$(2,16)$&$(8,8)$&$(4,2)$&$128$&$128,256,512,1024$\\
			$(2,67,53,13)$&$(2,4)$&$(2,8)$&$(2,2)$&$32$&$32,64$\\
			$(2,19,53,13)$&$(2,8)$&$(4,8)$&$(3,2)$&$64$&$64,128,256$\\
			\bottomrule
		\end{tabular}
	\end{table}
	
	These are necessary candidates for the arithmetic groups. Abstract sharpness does not prove realization of each candidate by these fields. For $(1,491,61,13)$, $m=4$, $t=2$ and the distinguished module has exponent $8<16$. Corollary~\ref{cor:typeII} therefore excludes family~II for both $G_1$ and $G$.
	
	We independently checked the finite quadratic class-number inputs by enumerating cycles of reduced ideals. For a positive fundamental discriminant $D$, represent a primitive reduced ideal by a pair $(a,b)$ with
	\begin{equation}\label{eq:reducedpairs} 0<a<\sqrt D,\qquad |\sqrt D-2a|<b<\sqrt D,\qquad b^2\equiv D\pmod{4a},\qquad \gcd\left(a,b,\frac{b^2-D}{4a}\right)=1. \end{equation}
	The reduction permutation is
	\begin{equation}\label{eq:reduction} \ell=\left\lfloor\frac{b+\sqrt D}{2a}\right\rfloor,\qquad b_1=2a\ell-b,\qquad (a,b)\longmapsto\left(\frac{D-b_1^2}{4a},b_1\right). \end{equation}
	Its cycles correspond to ordinary ideal classes; see \cite[Chapter~5]{Cohen1993}. Enumerating all pairs and cycles gives Table~\ref{tab:verification}. The bounds in \eqref{eq:reducedpairs} can be tested by integer square comparisons, and the floor in \eqref{eq:reduction} equals $\lfloor(b+\lfloor\sqrt D\rfloor)/(2a)\rfloor$, so the enumeration uses exact integer arithmetic. For $\Q(\sqrt{689})$ and $\Q(\sqrt{793})$ the same calculation gives ordinary class number $4$ in each case; their fundamental units are $105+4\sqrt{689}$ and $4393+156\sqrt{793}$, both of norm $1$. This verifies $t=2$ throughout Table~\ref{tab:examples} and checks the ordinary convention separately from the inherited biquadratic module structures. It does not replace the stabilization arguments for those modules.
	
	\begin{table}[ht]
		\centering\small
		\caption{Independent reduced-ideal counts for the quadratic inputs.}\label{tab:verification}
		\begin{tabular}{@{}lrrrr@{}}
			\toprule
			$(\eta,q,r,s)$&$D(F)$&reduced ideals&$h(F)$&$h_2(\eta qrs)$\\
			\midrule
			$(1,131,53,13)$&$361036$&$592$&$8$&$8$\\
			$(1,367,61,13)$&$1164124$&$864$&$16$&$16$\\
			$(1,491,61,13)$&$1557452$&$324$&$32$&$32$\\
			$(2,67,53,13)$&$369304$&$552$&$24$&$8$\\
			$(2,19,53,13)$&$104728$&$232$&$16$&$16$\\
			\bottomrule
		\end{tabular}
	\end{table}
	
	The quartic field needed for the first test can be specified directly from the radicals. Put
	\begin{equation} T=\Q(\sqrt{\eta q},\sqrt r,\sqrt s),\qquad T_\oo=TF_\oo. \end{equation}
	Since $K_\oo,K'_\oo,K''_\oo$ are the three distinct quadratic unramified extensions of $F_\oo$, their compositum is $T_\oo$ and has degree four over $F_\oo$. It lies in $\cH(F_\oo)$ and corresponds to the Frattini subgroup $\langle a^2,b^2\rangle$. Theorem~\ref{thm:arithmetic} gives
	\begin{equation} |X(T_\oo)|= \begin{cases} 2^{m+t-1},&|G|=2^{m+t+1},\\ 2^{m+t},&|G|>2^{m+t+1}. \end{cases} \end{equation}
	Thus an explicit triquadratic Iwasawa module decides whether the canonical first-layer quotient is the full tower; no choice of a metacyclic presentation enters the test.

	\subsection{Computations and group identification}
	The quartic test can be completed for every row of Table~\ref{tab:examples}. Put $T_0=T$ and $T_1=T(\sqrt2)$. Here $h(M)$ denotes the ordinary class number of a number field $M$.
	
	\begin{proposition}\label{prop:resolvedexamples} For the five tuples in Table~\ref{tab:resolvedexamples}, the ordinary class numbers of $T_0$ and $T_1$ are as listed, and \[ X(T_\oo)\cong C_4\times C_{2^{m-1}},\qquad [\cL(F_\oo):F_\oo]=2^{m+3}. \] In each case the full group has the presentation \[ G(F_\oo)\cong G(F)\cong\langle a,b\mid a^8=b^{2^m}=1,\ b^{-1}ab=a^{-1}\rangle. \] Moreover, \[ \cL(F)=\cH(K),\qquad \cL(F_\oo)=\cH(K)F_\oo=\cH(K_\oo). \] Thus the canonical first-layer quotient is the full tower group for all five examples. \end{proposition}
	
	\begin{table}[ht]
		\centering\small
		\caption{Completed arithmetic recovery from the quartic modules.}\label{tab:resolvedexamples}
		\begin{tabular}{@{}lrrrrr@{}}
			\toprule
			$(\eta,q,r,s)$&$m$&$h(T_0)$&$h(T_1)$&$|X(T_\oo)|$&$|G(F_\oo)|$\\
			\midrule
			$(1,131,53,13)$&$2$&$8$&$8$&$8$&$32$\\
			$(1,367,61,13)$&$3$&$48$&$48$&$16$&$64$\\
			$(1,491,61,13)$&$4$&$32$&$160$&$32$&$128$\\
			$(2,67,53,13)$&$2$&$24$&$24$&$8$&$32$\\
			$(2,19,53,13)$&$3$&$16$&$16$&$16$&$64$\\
			\bottomrule
		\end{tabular}
	\end{table}
	
	\begin{proof} The class numbers are computed using the multiquadratic class-number formula and exact unit-index calculations. In particular, $h_2(T_0)=h_2(T_1)=2^{m+1}$ in every row. In the third row the ordinary class number increases by a factor of $5$, while its $2$-part remains unchanged. \par We first define the cyclotomic layers and verify the ramification required for stabilization. For $\nu\ge0$, put $\Q_\nu=\Q(\zeta_{2^{\nu+2}}+\zeta_{2^{\nu+2}}^{-1})$, $T_\nu=T\Q_\nu$ and $F_\nu=F\Q_\nu$. The independent square classes of the defining radicands show that $\sqrt2\notin T$. Since every nontrivial finite subextension of $\Q_\oo/\Q$ contains $\Q(\sqrt2)$, we have $T\cap\Q_\oo=\Q$. Thus $[T_\nu:T]=[F_\nu:F]=2^\nu$, with $T_0=T$ and $T_1=T(\sqrt2)$. \par The square classes of $r$ and $s$ in $\Q_2^\times/\Q_2^{\times2}$ both equal that of $5$. If $\eta=1$, the square class of $q$ is that of $-5$ or $-1$, according as $q\equiv3$ or $7\pmod8$. If $\eta=2$, the square class of $2q$ is that of $-10$ or $-2$, respectively. Consequently every completion of $T_0$ above $2$ is isomorphic to $\Q_2(i,\sqrt5)$ when $\eta=1$, and to $\Q_2(\sqrt{-2},\sqrt5)$ when $\eta=2$. Adjoining $\sqrt2$ gives $\Q_2(\zeta_8,\sqrt5)$ in either case. For every $\nu\ge1$, the corresponding completion of $T_\nu$ is therefore $\Q_2(\zeta_{2^{\nu+2}},\sqrt5)$: adjoining $i$ to the real cyclotomic subfield gives the full local cyclotomic field. The extension $\Q_2(\sqrt5)/\Q_2$ is unramified, while $\Q_2(\zeta_{2^{\nu+2}})/\Q_2$ is totally ramified of degree $2^{\nu+1}$. Since the completion of $T_0$ has ramification index $2$ over $\Q_2$, the relative ramification index in $T_\nu/T_0$ is $2^\nu=[T_\nu:T_0]$. Hence every prime above $2$ is totally ramified from level zero. Only primes above $2$ can ramify in this cyclotomic tower. Proposition~\ref{prop:stability} now gives $|X(T_\oo)|=2^{m+1}$. Since $|X(K_\oo)|=2^{m+2}$, Theorem~\ref{thm:arithmetic} gives $k=1$, $|G(F_\oo)|=2^{m+3}$ and $\cL(F_\oo)=\cH(K_\oo)$. \par We next compare the finite and cyclotomic bases. Writing $D(M)$ for the absolute field discriminant, the quadratic-subfield discriminant product gives $D(T)=(4\eta qrs)^4=D(F)^4$. Since $[T:F]=4$, the relative discriminant of $T/F$ is trivial. Both fields are totally real, so $T/F$ is unramified at every finite and real place. Its base change $T_\nu/F_\nu$ is also unramified. Multiplicativity of ramification indices, together with total ramification in $T_\nu/T$, therefore gives ramification index $2^\nu$ in $F_\nu/F$ at every prime above $2$. Thus every nontrivial finite subextension of $F_\oo/F$ is ramified above $2$, whereas every finite subextension of $\cL(F)/F$ is unramified there. Consequently $F_\oo\cap\cL(F)=F$. \par The finite-level computations also certify $h_2(K)=2^{m+2}$. Since $K\subseteq T$, the extension $K/F$ is unramified, and hence $\cH(K)\subseteq\cL(F)$. Moreover, $[\cH(K):F]=2h_2(K)=2^{m+3}$. Disjointness from $F_\oo$ and unramified base change give $\cH(K)F_\oo\subseteq\cL(F_\oo)$ and $[\cH(K)F_\oo:F_\oo]=2^{m+3}$. Comparison with $|G(F_\oo)|$ proves $\cL(F_\oo)=\cH(K)F_\oo$. For every finite subextension $M/F$ of $\cL(F)/F$, the same disjointness and base-change argument gives $[M:F]=[MF_\oo:F_\oo]\le2^{m+3}$. Thus $\cL(F)/F$ is finite of degree at most $2^{m+3}$. The inclusion $\cH(K)\subseteq\cL(F)$ gives the reverse inequality, so $\cL(F)=\cH(K)$. We have therefore proved $\cL(F_\oo)=\cL(F)F_\oo=\cH(K)F_\oo=\cH(K_\oo)$. Restriction along this base change gives the isomorphism $G(F_\oo)\cong G(F)$. \par Finally, the ideal-extension calculations give \[ |\ker(\Cl_2(F)\longrightarrow\Cl_2(K))|=4 \] in every row. Under the preceding restriction isomorphism, the subgroup corresponding to $K_\oo=KF_\oo$ is identified with the subgroup corresponding to $K$. Transfer naturality and the transfer--capitulation correspondence therefore identify this finite capitulation kernel with the distinguished transfer kernel of $G(F_\oo)$. Proposition~\ref{prop:capitulation} selects family~I. Since $t=2$ and $k=1$, its parameters satisfy $\alpha=t+1=3$ and $n=m$, giving the asserted presentation. The subgroup corresponding to $T_\oo$ is then $\langle a^2,b^2\rangle\cong C_4\times C_{2^{m-1}}$. This subgroup is abelian, so the class-field correspondence gives $X(T_\oo)\cong C_4\times C_{2^{m-1}}$, as required. \end{proof}
	
	The class-number and capitulation calculations were performed in PARI/GP~2.15.4 \cite{PARI}, using exact unit-index computations and ideal-class extension maps.

	For the first row, the computations give $h(T)=h(T(\sqrt2))=8$. Stabilization therefore gives $|X(T_\oo)|=8$, and the quartic test recovers $[\cL(F_\oo):F_\oo]=32$. The distinguished capitulation kernel has order four, identifying family~I with $m=2$. Thus the quartic computation resolves the two candidates supplied by the first-layer data.
	
	For larger $m$, a quartic module distinguishes $k=1$ from $k\ge2$. At degree $2^j$, any of the three intermediate fields may then be used to determine whether the plateau continues. A module at degree $2^m$ always yields the full tower degree through
	\[ [\cL(F_\oo):F_\oo]=2^m|X(E)|. \]
	Although one such field has cyclic $X(E)$, the degree formula uses only its order and holds for all three choices.
	
	The computations above resolve all five rows of Table~\ref{tab:examples} using the first two cyclotomic levels. In each case the quartic module gives $k=1$, so no higher unramified module is needed to recover the tower degree. For other fields, finite-level computations of higher modules require class-number certification and a stabilization argument before Theorem~\ref{thm:main} can be applied.
	
	The remaining distinction is therefore precise. First-layer transfer data determine the canonical quotient in this class, but can leave every degree in the stated range possible. Higher abelian modules recover the degree, and at the upper endpoint they recover the full presentation. Away from that endpoint, equal orders need not determine an isomorphism class, as Proposition~\ref{prop:sharpness}(3) shows. The canonical filtration separates these questions and identifies the arithmetic data required for each.
	
	\section*{Acknowledgments}
	ChatGPT was used during the preparation of this article to assist with language editing and typographic errors. Responsibility for the mathematical content, calculations, and presentation remains with the author.

\end{document}